\documentclass[a4paper,11pt,leqno,twoside,final]{amsart}

\usepackage[latin1]{inputenc}               
\usepackage[OT1]{fontenc}          
\usepackage[usenames]{color}

\usepackage{array}                          
\usepackage{tabularx}
\usepackage[usenames]{colortbl}

\usepackage{amsmath,amssymb,amsthm}         
\usepackage[mathcal]{eucal}

\usepackage{verbatim,ae,ifthen,hhline,showkeys}                

\usepackage{a4wide}

\makeatletter
\renewcommand{\section}{\@startsection               
{section}{1}{0mm}{-\baselineskip}{0.5\baselineskip}{\raggedright\Large\bfseries\sffamily}}%
\makeatother

\newcommand{\caso}[2]{\textbf{Case #1:} #2.} 

\newcommand{\aut}{\textsc}
\newcommand{\whe}{\textsl}
\newcommand{\num}{\textbf}

\numberwithin{equation}{section}

\theoremstyle{plain}

\newtheoremstyle{slbody}%
{3pt}
{3pt}
{\normalfont\slshape}
{}
{\normalfont\bfseries\sffamily}
{.}
{ }
{}%

\newtheoremstyle{nmbody}%
{3pt}%
{3pt}%
{\normalfont}%
{}%
{\normalfont\bfseries\sffamily}%
{.}%
{ }%
{}%

\theoremstyle{slbody}

\newtheorem{theo}{Theorem}[section]
\newtheorem{lem}{Lemma}[section]
\newtheorem{cor}{Corollary}[section]

\theoremstyle{nmbody}

\renewcommand{\epsilon}{\varepsilon} 
\renewcommand{\theta}{\vartheta}
\renewcommand{\phi}{\varphi}

\newcommand{\ns}{\mathbf} 
\newcommand{\cs}{\mathrm} 

\newcommand{\supp}{\mathop{\rm supp}\nolimits}     

\newcommand{\np}{{\nabla_{+}}}     
\newcommand{\nm}{{\nabla_{-}}}
\newcommand{\npm}{{\nabla_{\pm}}}
\newcommand{\tp}{\mathop{\tau_{+}}\nolimits}       
\newcommand{\tm}{\mathop{\tau_{-}}\nolimits}
\newcommand{\tpm}{\mathop{\tau_{\pm}}\nolimits}

\newcommand{\de}{\partial}      
\newcommand{\di}{\, d}                             

\newcommand{\jp}[1]{\langle #1 \rangle}
\newcommand{\st}{\, : \,}                      

\newcommand{\q}{\left} 
\newcommand{\p}{\right}

\renewcommand{\leq}{\leqslant} 
\renewcommand{\geq}{\geqslant}

\newcommand{\prf}{\textsc{Proof.} } 
\newcommand{\bs}{~$\blacksquare$}

\newcommand{\beq}{\begin{equation}}
\newcommand{\eeq}{\end{equation}}

\begin{document}

\title[A dispersive estimate for the linear wave equation]%
{A dispersive estimate for the linear wave equation with an
electromagnetic potential}

\author{Davide Catania}
   \address{ Davide Catania,
Dipartimento di Matematica, Universit\`a di Pisa, Largo B.
Pontecorvo 5, 56127 Pisa, Italy;  e--mail: {\normalfont
catania@mail.dm.unipi.it}}


\begin{abstract}
We consider radial solutions to the Cauchy problem for the linear
wave equation with a small short--range electromagnetic potential
(the ``square version'' of the massless Dirac equation with a
potential) and zero initial data. We prove two a priori estimates
that imply, in particular, a dispersive estimate.
\end{abstract}

\maketitle

\section{Introduction}

\bigskip

In this paper, we investigate the dispersive properties of the
linear wave equation with an electromagnetic potential, that is
 \beq
 \square_A u  = F \qquad \qquad (t,x) \in \ns [0,\infty[\times \ns R^3,
 \eeq
where
 \begin{align}
 & x= (x_1, x_2,x_3), \qquad \qquad r = |x|,\\
 & \square_A  = \square - A\cdot \nabla_{t,r}, \\
 & \square = \de_t^2-\Delta = \de_t^2-(\de_{x_1}^2+\de_{x_2}^2+\de_{x_3}^2), \\
 & \nabla_{t,r} = \begin{pmatrix} \de_t \\
 \de_r \end{pmatrix}.
 \end{align}
The fact that the potential $A=A(t,x)$ is electromagnetic means that
$A \in i\ns R \times i \ns R$, where $i$ is the imaginary unit. This
will play a crucial role in the development of the proof, since
electromagnetic potential are gauge invariant (see what follows).

We restrict ourselves to radial solutions $u=u(t,r)$, with
$F=F(t,r)$ and
 \beq
 A=A(t,r) = \binom{A_0(t,r)}{A_1(t,r)}, \qquad A_0, A_1 \in i\ns
 R. \eeq
We assume further that the potential decreases sufficiently rapidly
when $r$ approaches infinity; more precisely, we suppose that
 \beq \label{cond.A}
 \sum_{j \in \ns Z} 2^{-j}\jp{2^{-j}}^{\epsilon_A} ||\phi_j A||_{L^\infty_r}  \leq
 \delta_A
 \eeq
(that is, $A$ is a short--range potential), where $\epsilon_A>0$,
$\delta_A$ is a sufficiently small positive constant independent of
$r$ (see Section \ref{sec.apriori})  and the sequence $(\phi_j)_{j
\in \ns Z}$ is a Paley--Littlewood partition of unity, which means
that $\phi_j(r) = \phi (2^j r)$ and  $\phi : \ns R^+ \longrightarrow
\ns R^+$ ($\ns R^+$ is the set of all non--negative real numbers) is
a function so that
 \begin{enumerate}
 \item $\supp \phi =\{ r \in \ns R \st 2^{-1} \leq r \leq 2 \}$;
 \item $\phi(r)>0 \quad $ for $\quad 2^{-1} < r < 2$;
 \item $\sum_{j \in \ns Z} \phi (2^j r) = 1 \quad $ for each $\quad r \in \ns R^+$.
 \end{enumerate}
In other words, $\sum_{j \in \ns Z} \phi_j (r) = 1$ for all $r \in
\ns R^+$ and
 \beq
 \supp \phi_j = \{ r \in \ns R \st 2^{-j-1} \leq r \leq
2^{-j+1} \}.
 \eeq

\medskip

It is well--known that there exists a unique global solution to the
Cauchy problem
 \beq \label{prob.base} \begin{cases}
 \square_A u  = F & \qquad (t,x) \in \ns [0,\infty[\times \ns R^3, \\
 u(0,x)=\de_t u(0,x)=0 & \qquad x \in \ns R^3;
 \end{cases}
 \eeq
in particular, this fact holds for the smaller class of radial
solutions, that is for the problem
 \beq \label{prob.rad} \begin{cases}
 \square_A u  = F & \qquad (t,r) \in \ns [0,\infty[\times \ns R^+, \\
 u(0,r)=\de_t u(0,r)=0 & \qquad r \in \ns R^+.
 \end{cases}
 \eeq

Let introduce the change of coordinates
 \beq
 \tpm{} \doteq \frac{t\pm r}{2}
 \eeq
and the standard notation $\jp{s} \doteq \sqrt{1+s^2}$; our main
result can be expressed as follows.

\begin{theo} \label{teo.est}
Let $u$ be a radial solution to \eqref{prob.base}, i.e. a solution
to \eqref{prob.rad}, where $A=A(t,r)$ is an electromagnetic
potential satisfying \eqref{cond.A} for some $\delta_A >0$ and
$\epsilon_A>0$. Then, for every $\epsilon>0$, there exist two
positive constants $\delta$ and $C$ (depending on $\epsilon$) such
that for each $\delta_A \in ]0,\delta]$, one has
 \beq
 || \tp  u ||_{\cs L^\infty_{t,r}} \leq C ||
 {\tp}r^2\jp{r}^{\epsilon}
 F ||_{\cs L^\infty_{t,r}}.
 \eeq
\end{theo}

Let introduce the differential operators
 \beq
 \npm \doteq  \de_t \pm \de_r.
 \eeq
The proof of the previous a priori estimate follows easily from the
following one.

\begin{lem} \label{lem.est} Under the same conditions of Theorem \ref{teo.est}, for every $\epsilon>0$, there exist two
positive constants $\delta$ and $C$ (depending on $\epsilon$) such
that for each $\delta_A \in ]0,\delta]$, one has
 \beq \label{eq.nmufin}
 || \tp r \nm u ||_{\cs L^\infty_{t,r}} \leq C ||
 {\tp}r^2\jp{r}^{\epsilon}
 F ||_{\cs L^\infty_{t,r}}.
 \eeq
\end{lem}

An immediate consequence of Theorem \ref{teo.est} is the following
dispersive estimate.

\begin{cor} \label{cor.disp} Under the same conditions of Theorem \ref{teo.est},
 for every $\epsilon>0$, there exist two
positive constants $\delta$ and $C$ (depending on $\epsilon$) such
that for each $\delta_A \in ]0,\delta]$, one has
 \beq
  |u(t,x)| \leq \frac{C}{ t}   || \tp r^2 \jp{r}^{\epsilon} F ||_{\cs
 L^\infty_{t,r}}
 \eeq
for every $t>0$.
\end{cor}

The idea to prove the lemma is the following. First of all, the
potential term in \eqref{prob.rad} can be thought as part of the
forcing term, that is $\square_A u = F$ can be viewed as
 \beq
 \square u = F_1 \doteq F + A\cdot \nabla_{t,r} u.
 \eeq
Then we can rewrite this equation in terms of $\tpm{}$ and $\npm$
(see Section \ref{sec.apriori}), obtaining
 \beq
 \np\nm v = G,
 \eeq
where
 \beq \label{def.vG}
 v(t,r) \doteq r u(t,r) \qquad \text{and} \qquad G(t,r) \doteq r
 F_1(t,r).
 \eeq
This last equation can be easily integrated to obtain a relatively
simple explicit representation of $(\nm v)(\tp,\tm)$ in terms of
$G$.

Another fundamental step consists in taking advantage of the gauge
invariance property of the electromagnetic potential $A$, which
means that, set
 \beq
 A_\pm \doteq \frac{A_0 \pm A_1}{2},
 \eeq
we can assume, with no loss of generality, that $A_+ \equiv 0$ (see
\cite{bs}, p. 34). This implies that
 \beq
 A \cdot \nabla_{t,r} u = A_- \nm u + A_+ \np u = A_- \nm u
 \eeq
and hence
 \beq \label{F1Am}
 F_1 = F+ A_- \nm u ,
 \eeq
thus
 \beq \label{eq.G-}
 G = rF + r A_- \nm u  = rF + A_- \nm v + \frac{1}{r} A_- v.
 \eeq
Obviously, one has
  \beq \label{cond.Am}
 \sum_{j \in \ns Z} 2^{-j}\jp{2^{-j}}^{\epsilon_A} ||\phi_j A_-||_{L^\infty_r}  \leq
 \delta_A.
 \eeq
These simplifications, combined with technical Lemma \ref{l.1} and
the estimate of Lemma \ref{l.2}, allow us to easily obtain Lemma
\ref{lem.est} and Theorem \ref{teo.est}.

\medskip

The dispersive properties of evolution equations are important for
their physical meaning and, consequently, they have been deeply
studied, though the problem in its generality is still open. The
dispersive estimate obtained in Corollary \ref{cor.disp} provides
the natural decay rate, that is the same rate one has for the
non--perturbed wave equation (see \cite{give,kt}), i.e.  a
$t^{-(n-1)/2}$ decay in time, where $n$ is the space dimension (in
our case, $n=3$). The generalization to the case of a
potential--like perturbation has been considered widely (see
\cite{beals,bpstz,cuccagna,cs,ghk,gv,pstz,visciglia1,visciglia2,
yajima95,yajima99}), also for the Schr\"odinger equation (see
\cite{gst,gt,jss,rs}). Recently, D'Ancona and Fanelli have
considered in \cite{daf} the case
 \beq \begin{cases}
 \de_t^2 u(t,x) + H u =0, & \quad (t,x) \in
 \ns R \times \ns R^3, \\
 u(0,x)=0, \quad \de_t u(0,x)=g(x),
 \end{cases} \eeq
where
 \begin{gather}
 H \doteq - (\nabla+iA(x))^2 + B(x), \\
 A : \ns R^3 \longrightarrow \ns R^3, \qquad B : \ns R^3 \longrightarrow \ns
 R.
 \end{gather}
Under suitable condition on $A$, $\nabla A$ and $B$, in particular
 \beq \label{cond.daf}
 |A(x)| \leq \frac{C_0}{r\jp{r}(1+|\lg r|)^\beta},
 \qquad \sum_{j=1}^3 |\de_j A_j (x)| + |B(x)| \leq \frac{C_0}{r^2(1+|\lg r|)^\beta},
 \eeq
with $C_0>0$ sufficiently small, $\beta >1$ and $r=|x|$, they have
obtained the dispersive estimate
 \beq
 |u(t,x)| \leq \frac{C}{t} \sum_{j \geq 0} 2^{2j} || \jp{r}
 w_\beta^{1/2} \phi_j(\sqrt{H}) g ||_{L^2},
 \eeq
where $w_\beta \doteq r(1+|\log r|)^\beta$ and $(\phi_j)_{j\geq 0}$
is a non--homogeneous Paley--Littlewood partition of unity on $\ns
R^3$.

In this paper, restricting ourselves to radial solutions, we are
able to obtain the result in Corollary \ref{cor.disp}, which is
optimal from the point of view of the estimate decay rate $t^{-1}$
and improve essentially the assumptions on the potential, assuming
weaker condition \eqref{cond.A} instead of \eqref{cond.daf}.

\section{A priori estimates} 
\label{sec.apriori}

First of all, we reformulate our problem  taking advantage of the
radiality of the solution $u$ to \eqref{prob.rad}. Indeed, since
$\Delta_{\ns S^2} u(t,r) =0$ and $v=ru$, we have
\begin{align}
 \square u(t,r) & = ( \de_t^2 - \Delta_x) u = \q(\de_t^2 - \de_r^2 -
\frac{2}{r}\de_r -\frac{1}{r^2}\Delta_{\ns S^2}\p) u(t,r) \\
& = \frac{1}{r} \de_t^2 v(t,r) - \frac{1}{r} \de_r^2 v(t,r) \\
& = \frac{1}{r} \np\nm v(t,r) = \frac{1}{r}\nm\np v(t,r).
\end{align}
Recalling \eqref{def.vG} and \eqref{F1Am}, we get that the equation
in \eqref{prob.rad} is equivalent to
  \beq
 \np\nm v = G.
 \eeq
Let us notice that the support of $u(t,r)$ is contained in the
domain $\{ (t,r) \in \ns R^2 : r>0, \, t>r \}$, therefore we have
 \beq \label{eq.suppv} \supp v(\tp,\tm) \subseteq \{ (\tp,\tm) \in \ns R^2 : \tm>0, \, \tp > \tm \}. \eeq

From this fact, we get
 \beq
 \nm v(\tp,\tm) = \nm v(\tm,\tm) + \int_{\tm}^{\tp} G(s,\tm) \di
 s = \int_{\tm}^{\tp} G(s,\tm) \di
 s.
 \eeq
Let us observe that,  for each $s \in [\tm,\tp]$, we have
 \beq
 s \leq \tp, \qquad \qquad {s-\tm} \leq \tp - \tm = r,
 \eeq
hence
 \begin{align*}
 \q| \int_{\tm}^{\tp} G(s,\tm) \di s \p| & \leq  \int_{\tm}^{\tp}
 \frac{s\jp{s-\tm}^\epsilon |G(s,\tm)|}{\jp{s}\jp{s-\tm}^\epsilon} \di
 s  \\
 & \leq || \tp\jp{r}^\epsilon G ||_{\cs L^\infty_{t,r}}
 \int_{\tm}^{\tp} \jp{s}^{-1} \jp{s-\tm}^{-\epsilon} \di s
 \end{align*}
for every $\epsilon > 0$. Applying lemma \ref{l.1} (see the end of
this section), we conclude
 \beq
 \tp | \nm v(\tp,\tm) | \leq Cr || \tp \jp{r}^\epsilon G
 ||_{\cs L^\infty_{t,r}};
 \eeq
recalling that $G$ satisfies \eqref{eq.G-},
  we obtain
 \begin{equation} \label{eq.tpnmv} \begin{split}
 \tp | \nm v(\tp,\tm) | & \leq C r \Big( || \tp \jp{r}^\epsilon  A_-
 \nm v  ||_{\cs L^\infty_{t,r}} + || \tp \jp{r}^\epsilon r^{-1}  A_-
 v  ||_{\cs L^\infty_{t,r}}  \\
 & \hspace{1cm}   +  || \tp \jp{r}^\epsilon rF
 ||_{\cs L^\infty_{t,r}} \Big).
 \end{split} \end{equation}

Now,  if we choose for the moment $\epsilon \leq \epsilon_A$, we
have
 \beq
  r\jp{r}^\epsilon \phi_j(r) |A_-(t,r)| \leq C
   2^{-j}\jp{2^{-j}}^{\epsilon_A} || \phi_j A_- ||_{\cs L^\infty_r}
  \eeq
(here and in the following, we assume that $C=C(\epsilon)>0$ could
change time by time), thus
 \beq \label{eq.nmvI} \begin{split}
 r|| \tp \jp{r}^\epsilon A_-
 \nm v  ||_{\cs L^\infty_{t,r}} & \leq C || \tp \nm v
 ||_{\cs L^\infty_{t,r}} \sum_{j \in \ns Z} 2^{-j}\jp{2^{-j}}^{\epsilon_A}
  || \phi_j A_- ||_{\cs L^\infty_r} \\
  & \leq C\delta_A || \tp \nm v
 ||_{\cs L^\infty_{t,r}}, \end{split}
 \eeq
where we have used the fact that $(\phi_j)_{j \in \ns Z}$ is a
Paley--Littlewood partition of unity and property \eqref{cond.Am}.

Moreover, $v(\tp,\tp)=0$ because of \eqref{eq.suppv}, whence
 \beq
 v(\tp,\tm) = -\int_{\tm}^{\tp} \nm v (\tp,s) \di s
 \eeq
and consequently
 \beq \label{est.vnv}
 |v(\tp,\tm)| \leq \int_{\tm}^{\tp} | \nm v (\tp,s) | \di s \leq r ||
 \nm v  ||_{\cs L^\infty_{t,r}}.
 \eeq
Thus we have
 \beq
 \jp{r}^\epsilon \phi_j(r)|A_-(t,r)v(\tp,\tm)| \leq C 2^{-j}\jp{2^{-j}}^{\epsilon_A}
  || \phi_j A_- ||_{\cs L^\infty_r} ||
 \nm v  ||_{\cs L^\infty_{t,r}},
 \eeq
which implies
 \beq \label{eq.nmvII}
 r|| \tp \jp{r}^\epsilon r^{-1} A_-
 v  ||_{\cs L^\infty_{t,r}} \leq C\delta_A || \tp \nm v
 ||_{\cs L^\infty_{t,r}}.
 \eeq
Using \eqref{eq.nmvI} and \eqref{eq.nmvII} in \eqref{eq.tpnmv}, we
deduce
 \beq \label{eq.nmvfin}
 || \tp  \nm v ||_{\cs L^\infty_{t,r}} \leq C || \tp r^2\jp{r}^\epsilon F
 ||_{\cs L^\infty_{t,r}},
 \eeq
provided $\delta_A$ is sufficiently small. For instance, one can
take $\delta_A$ such that $3C^2\delta_A\leq 1$.

From the definition of $v$, we have
 \beq \label{eq.nmuv}
 r\nm u = \nm v + u
 \eeq
and hence
 \beq
 | \tp r \nm u | \leq | \tp \nm v | + | \tp u |,
 \eeq
so
 \begin{align*}
 || \tp r \nm u ||_{\cs L^\infty_{t,r}} & \leq || \tp \nm v ||_{\cs
 L^\infty_{t,r}} + || \tp u ||_{\cs L^\infty_{t,r}} \\
 & \leq C \q( || \tp r^2\jp{r}^\epsilon F  ||_{\cs L^\infty_{t,r}}
 + || \tp r^2\jp{r}^\epsilon F_1  ||_{\cs L^\infty_{t,r}} \p),
 \end{align*}
where we have used \eqref{eq.nmvfin} and the inequality in Lemma
\ref{l.2}. But
 \beq \begin{split}
 r^2\jp{r}^\epsilon |F_1| & \leq r^2\jp{r}^\epsilon |A_-|\cdot |\nm u| +
 r^2\jp{r}^\epsilon |F| \\
 & \leq \q( \sum_{j \in \ns Z} r\jp{r}^{\epsilon_A} \phi_j |A_-|
 \p) r |\nm u|
  +  r^2\jp{r}^\epsilon |F| \\
 & \leq C\delta_A r |\nm u|
  +  r^2\jp{r}^\epsilon |F|,
 \end{split} \eeq
thus
 \beq
 || \tp r^2\jp{r}^\epsilon F_1  ||_{\cs L^\infty_{t,r}} \leq
 C\delta_A || \tp r \nm u ||_{\cs L^\infty_{t,r}} + || \tp r^2\jp{r}^\epsilon
 F ||_{\cs L^\infty_{t,r}}
 \eeq
and consequently
 \beq
 || \tp r \nm u ||_{\cs L^\infty_{t,r}} \leq C ||
 {\tp}r^2\jp{r}^{\epsilon}
 F ||_{\cs L^\infty_{t,r}}
 \eeq
provided $\delta_A>0$ small enough, that is Lemma \ref{lem.est}.

Now we use the fact that, because of \eqref{eq.nmuv}, we have
 \beq
 |\tp u| \leq | \tp r \nm u| + | \tp \nm v |;
 \eeq
combining this estimate with \eqref{eq.nmufin} and
\eqref{eq.nmvfin}, we finally conclude
 \beq
 || \tp u ||_{\cs L^\infty_{t,r}} \leq C ||
 {\tp}r^2\jp{r}^{\epsilon}
 F ||_{\cs L^\infty_{t,r}},
 \eeq
and also Theorem \ref{teo.est} is proven.

\begin{lem} \label{l.1} For each $\epsilon >0$, there exists a positive constant
$C=C(\epsilon)$ such that
 $$
 \int_{\tm}^{\tp} \jp{s}^{-1} \jp{s-\tm}^{-\epsilon} \di s \leq
 \frac{Cr}{\tp}
 $$
\end{lem}

\prf We distinguish two cases.

\caso{1}{$\tp \geq 2\tm$} Let us notice that since $r = \tp - \tm
\geq \tp/2$, in this case it is sufficient to prove that
 \beq
 \int_{\tm}^{\tp} \jp{s}^{-1} \jp{s-\tm}^{-\epsilon} \di s \leq C_0(\epsilon).
 \eeq
We observe that $s-\tm \geq s/2$ provided $s \geq 2\tm$, so
 \begin{align*}
 \int_{\tm}^{\tp} \jp{s}^{-1} \jp{s-\tm}^{-\epsilon} \di s & \leq
 \int_{\tm}^{\tm +1} \jp{s}^{-1} \di s + 2^{\epsilon} \int_{\tm
 +1}^{\tp+1} s^{-(1+\epsilon)} \di s \\
 & \leq \frac{1}{\jp{\tm}} + 2^{\epsilon} \int_1^\infty
 s^{-(1+\epsilon)} \di s \\
 & \leq 1 + C_1(\epsilon).
 \end{align*}

\caso{2}{$\tp < 2\tm$} We use the estimates $\jp{s}^{-1} < 2/\tp$
and $\jp{s-\tm}^{-\epsilon} \leq 1$ to get
 \beq
 \int_{\tm}^{\tp} \jp{s}^{-1} \jp{s-\tm}^{-\epsilon} \di s \leq \frac{2}{\tp}(\tp
 - \tm) = \frac{2r}{\tp}. \eeq
This concludes the proof.\bs

\medskip

In the case $A\equiv 0$ (non--perturbed equation), we have the
following version of the estimate in Theorem \ref{teo.est}. It
consists in a slight modification of estimate $(1.8)$ shown in
\cite{ghk}, p. 2269.

\begin{lem} \label{l.2} Let $u$ be the solution to
 \beq \begin{cases}
 \square u  = F & \qquad (t,r) \in \ns [0,\infty[\times \ns R^+, \\
 u(0,r)=\de_t u(0,r)=0 & \qquad r \in \ns R^+.
 \end{cases}
 \eeq
Then, for every $\epsilon >0$, there exists $C>0$ such that \beq
  || \tp  u ||_{\cs L^\infty_{t,r}} \leq C || \tp r^2\jp{r}^\epsilon
  F
 ||_{\cs L^\infty_{t,r}}.
 \eeq
\end{lem}

\prf Let notice that $u$ is the solution to \eqref{prob.rad} with $A
\equiv 0$. Then, from \eqref{eq.tpnmv}, we have
 \begin{equation}
 \tp | \nm v(\tp,\tm) |  \leq C || \tp r^2 \jp{r}^\epsilon F
 ||_{\cs L^\infty_{t,r}},
 \end{equation}
where $v=ru$. Using \eqref{est.vnv}, we deduce
 \beq
 \tp |u| = \tp |v|r^{-1} \leq ||\tp \nm v ||_{\cs L^\infty_{t,r}}
 \eeq
and hence the claim.\bs

\bibliographystyle{amsalpha}

\end{document}